\documentclass[11pt]{amsart}

\usepackage{amsmath,amssymb,amscd,amsthm,amsxtra}

\usepackage{graphicx, mathrsfs}
\usepackage{amssymb,amsmath,amsthm}
\usepackage{mathrsfs,dsfont}

\usepackage{amsmath}
\usepackage{amsthm}
\usepackage{latexsym}
\usepackage{amssymb,amsfonts}
\usepackage{xcolor}
\usepackage{marvosym}
\usepackage{bbm}
\usepackage{graphicx, color} 
\usepackage{palatino, url, multicol}

\usepackage[all]{xy}

\makeatletter
\DeclareRobustCommand\widecheck[1]{{\mathpalette\@widecheck{#1}}}
\def\@widecheck#1#2{%
   \setbox\z@\hbox{\m@th$#1#2$}%
   \setbox\tw@\hbox{\m@th$#1%
      \widehat{%
         \vrule\@width\z@\@height\ht\z@
         \vrule\@height\z@\@width\wd\z@}$}%
   \dp\tw@-\ht\z@
   \@tempdima\ht\z@ \advance\@tempdima2\ht\tw@ \divide\@tempdima\thr@@
   \setbox\tw@\hbox{%
      \raise\@tempdima\hbox{\scalebox{1}[-1]{\lower\@tempdima\box\tw@}}}%
   {\ooalign{\box\tw@ \cr \box\z@}}}
\makeatother

\allowdisplaybreaks[2]

\newtheorem{theorem}{Theorem} [section]

\newtheorem{remark}[theorem]{Remark}

\newcommand{\Z}{\mathbb{Z}}

\newcommand{\R}{\mathbb{R}}
\newcommand{\C}{\mathbb{C}}
\newcommand{\T}{\mathbb{T}}
\newcommand{\E}{\mathbb{E}}

\newcommand{\EE}{\mathcal E}

\newcommand{\ft}{\widehat}

\newcommand{\cj}{\overline}

\numberwithin{equation}{section}

\begin{document}

\title[Gauge Transformations and Gaussian measures]{Absolute continuity of Brownian bridges under certain gauge transformations}
\author[Nahmod]{Andrea R. Nahmod$^1$}
\address{$^1$ Department of Mathematics \\ University of Massachusetts\\  710 N. Pleasant Street, Amherst MA 01003}
\email{nahmod@math.umass.edu}
\thanks{$^1$ The first author is funded in part by NSF DMS 0803160 and a 2009-2010 Radcliffe Institute for Advanced Study Fellowship.}

\author[Rey-Bellet]{Luc Rey-Bellet$^2$}
\address{$^2$ Department of Mathematics \\
University of Massachusetts\\ 710 N. Pleasant Street, Amherst MA 01003 }
\email{luc@math.umass.edu}
\thanks{$^2$ The  second author is funded in part by NSF DMS 0605058}

\author[Sheffield]{Scott Sheffield$^3$}
\thanks{$^3$ The third author is funded in part by NSF CAREER Award DMS 0645585 and a Presidential Early Career Award for Scientists and Engineers (PECASE)} 
\address{$^3$ Department of Mathematics\\
Massachusetts Institute of Technology\\ 77 Massachusetts Avenue,  Cambridge, MA 02139}
\email{sheffield@math.mit.edu}

\author[Staffilani]{Gigliola Staffilani$^4$}
\address{$^4$  Department of Mathematics\\
Massachusetts Institute of Technology\\ 77 Massachusetts Avenue,  Cambridge, MA 02139}
\email{gigliola@math.mit.edu}
\thanks{$^4$ The fourth author is funded in part by NSF
DMS 0602678 and a 2009-2010 Radcliffe Institute for Advance Study Fellowship.}
\date{}
\begin{abstract}
We prove absolute continuity of Gaussian measures associated to complex Brownian bridges
under certain gauge transformations. As an application we prove that the invariant measure for 
the periodic  derivative nonlinear Schr\"odinger equation  obtained  by  Nahmod, Oh, Rey-Bellet and 
Staffilani in \cite{NORBS}, and with respect to which they proved  almost surely global well-posedness,  
coincides with the weighted Wiener measure constructed by Thomann and 
Tzvetkov \cite{TTzv}. Thus, in  particular we prove the invariance of the measure constructed in \cite{TTzv}. 
\end{abstract}
\maketitle

\section{Introduction}\label{intro} 
This note is a continuation of the paper \cite{NORBS}. There we  constructed an invariant measure for the periodic
derivative nonlinear Schr\"odinger equation (DNLS) \eqref{DNLS} in one dimension and established global 
well-posedness, almost surely for data living in its support.  This was achieved by introducing a gauge 
transformation $G$, see \eqref{thegauge}, and by considering the {\it gauged} DNLS equation 
(GDNLS) \eqref{GDNLS+} in order to obtain the necessary  estimates. We constructed a
weighted Wiener measure $\mu$, which we proved to be invariant under the flow  of the GDNLS equation,  and used 
it to show the almost surely global well-posedness for  the GDNLS initial value problem, in particular almost 
surely for data in a certain Fourier-Lebesgue space  scaling like $H^{\frac{1}{2}-\epsilon}(\T),$ 
for  small  $\epsilon >0$.  To go back to the original DNLS equation we applied the inverse  
gauge transformation $G^{-1}$ and obtained an invariant measure  $\mu\circ G =: \gamma$  with respect 
to which almost surely global well-posedness is then proved for the  DNLS Cauchy initial value problem\footnote[5]{In \cite{NORBS} $\mu\circ G$ is called $\nu$; here we relabel it $\gamma$ to a priori distinguish it
from the name we give to the one constructed in \cite{TTzv}.}.
On the other hand,  Thomann and Tzvetkov \cite{TTzv}  constructed a weighted Wiener measure $\nu$ 
and proposed it as a natural candidate for an invariant measure for the DNLS equation. 
A natural question, left open in \cite{NORBS}, is the absolute continuity of  the two measures $\gamma$ and 
$\nu$  or  equivalently,  the absolutely continuity of $\mu$ and $\nu \circ G^{-1}$. 
As shown in this note, this question is easily answered after one understands the absolute continuity 
between Gaussian measures  naturally associated with complex Brownian bridges and their images 
under certain gauge transformations such as $G$. This is the heart of the matter of this note.  At the end 
we prove  that $\mu=\nu\circ G^{-1}$  (or equivalently that $\gamma = \nu$) thus in particular establishing 
the invariance of  the measure $\nu$ constructed in \cite{TTzv}, see Theorem \ref{maintheorem1}. 
Our results follow by combining the results on global well-posedness  and invariant measure for 
GDNLS  \eqref{GDNLS+} obtained by Nahmod, Oh, Rey-Bellet and Staffilani in  \cite{NORBS} with the 
explicit computation of the image of the measure under the gauge transformation. The  key to understand 
the latter is to actually understand how the Gaussian part  of the measure changes under the 
gauge since the transformation of the  weight is computed easily (see subsection \ref{heuristic} below). This is 
achieved in Theorem \ref{T1} in Section \ref{general} of this paper.

Certainly there is a vast literature on the topic of Gaussian measures under  nonlinear transformations 
\cite{CM,  Ramer, Kusuoka, Buck,  Enchev, Kallianpur, Cambronero} as well 
as \cite{Bogachev} and other references therein. But as we will show below the nature of the 
gauge transformation $G$ does not fit in the context  of these works and a different approach needs to be 
introduced.   For many nonlinear partial differential equations gauge transformations are an essential tool 
to convert one kind of  nonlinearity into another one, where resonant interactions are more manageable and 
hence estimates can be proved. Therefore we expect the general nature of  the central theorem of this note,
Theorem  \ref{T1},  as well as some of the ideas  behind its proof,  to be applicable in other situations beyond the 
DNLS context.
  
\section{Invariance of weighted Wiener measure for DNLS}\label{InvarianceDNLS} 
\smallskip
As stated in the introduction our motivation arises from the recent paper  by Nahmod, Oh, Rey-Bellet, and Staffilani \cite{NORBS} we recall now the set up of that paper and formulate the problem that we want to solve here in that context. We consider the derivative 
nonlinear  Schr\"odinger equation (DNLS) on the circle $\mathbb{T}$, i.e., 
\begin{equation} \label{DNLS}
\begin{cases} u_t (x,t)  \, - i \,  u_{xx}(x,t)  \, = \, \lambda \left(|u|^2(x,t) u(x,t)\right)_x \\
u(x,0) = u_0(x), 
\end{cases}
\end{equation}
where $x \in \mathbb{T}$, $t \in \mathbb{R}$ and $\lambda\in \R$ is fixed. In the rest of the paper we will set $\lambda=1$ for simplicity.    Our goal is to show that this problem
defines a dynamical system, in the sense of ergodic theory.  Let us denote by $\Psi(t)$  
the flow map associated to our nonlinear equation, i.e., the solution of \eqref{DNLS}, whenever it exists,  
is given by   $u(x,t) = \Psi(t)(u_0(x))$.    Let further  $(\mathcal{B}, \mathcal{F}, \nu)$ be a probability space
where $\mathcal{B}$ is a space (here $\mathcal{B}$ will be a separable Banach space),  $\mathcal{F}$ is a 
$\sigma$-algebra (here it will always be the Borel $\sigma$-algebra) and  $\nu$ is a probability measure.      
The flow map $\Psi(t)$ define a dynamical system on the probability space $(\mathcal{B}, \mathcal{F}, \nu)$ 
if 

\smallskip

\noindent{\bf (a)} ($\nu$-almost sure wellposedness.)  There exists a subset $\Omega \subset \mathcal{B}$ with 
$\nu(\Omega)=1$ such that the flow map $\Psi(t): \Omega \to \Omega$ is well defined and continuous  in $t$ 
for all $t \in \mathbb{R}$. 

\smallskip

\noindent{\bf (b)} (Invariance of the measure $\nu$.)  The measure $\nu$ is invariant under the flow $\Phi(t)$, 
i.e., 
$$ 
\int f\left( \Psi(t) (u)\right) \,  d \nu(u) \, =\, \int f(u)   d \nu(u) \,,
$$     
for all $f  \in L^1(\mathcal{B}, \mathcal{F}, \nu)$ and all $t \in \mathbb{R}$.

The measure $\nu$ here, in a sense,  is a substitute for a conserved quantity and in fact $\nu$ is 
constructed by using  a certain conserved quantity.  Recall that the DNLS equation \eqref{DNLS}  is 
completely integrable  (c.f. \cite{KaupNew, HayOz2})  and among the conserved quantities are 
\begin{eqnarray}
\mbox{Mass:}  \qquad {m(u)} &=& \frac{1}{2\pi} \int  |u|^2 \, dx \,, \label{mass}    \\
\mbox{ Energy:}  \qquad {E(u)} &=& \int | u_x|^2 \, dx + \frac{3}{2}  {\rm Im} \int  u^2 {\cj {u}\, \cj{u_x}} \, dx + 
\frac{1}{2} \int |u|^6 \, dx \,,  \label{energy} \\
\mbox{ Hamiltonian:} \qquad {H(u)} &=& {\rm Im} \int u {\cj u}_x \, dx  + \frac{1}{2} \int  |u|^4 \, dx\,.         
\label{Hamiltonian}
\end{eqnarray}

We consider a probability measure $\nu$ which is based on the conserved quantity $E(u)$ 
(as well as the mass $m(u)$).  Let us decompose $u= a+ ib$ into real and imaginary part, and let 
us consider first the purely formal but suggestive expression for $\nu$.
 
\begin{align} \label{GibbsUG}
 d \nu & = C^{-1}  \chi_{\{\|u\|_{L^2} \le B\}}  e^{-\frac{\beta}{2} N(u)}
e^{-\frac{\beta}{2} \int (|u|^2 + |u_x|^2) dx} \prod_{x \in\T} da(x) db(x).
\end{align}
where 
\begin{equation}\label{nlinenergy}
N(u)= \frac{3}{2}  {\rm Im} \int  u^2 {\cj {u}\, \cj{u_x}} \, dx +  \frac{1}{2} \int |u|^6 \, dx 
\end{equation}
is the non-quadratic part of the energy $E(u)$. 
Note that we have added the conserved quantity $\int |u|^2 dx$  to the quadratic part of $E(u)$ such as 
to make it positive definite. The constant $\beta > 0$ does not play any  particular role here and, for simplicity,  
we choose $\beta=1$.  Note however that all the measures for different $\beta$ would  be  invariant under the 
flow and they are all mutually singular.  The cutoff on the $L^2$-norm is necessary to make the measure
normalizable since $N(u)$ is not bounded below. We will also see that this measure is  well-defined only for 
$B$ under a certain critical value $B^*$. 

The expression in \eqref{GibbsUG} at this stage is purely formal since there is no Lebesgue measure in 
infinite  dimensions. In order to give a rigorous definition of the measure  $\nu$  in \eqref{GibbsUG} 
one needs to:

\smallskip
\noindent {\bf (a)} Make sense of the Gaussian part of the measure \eqref{GibbsUG}, that is  of 
the formal expression
\begin{align} \label{WienerUG}
d \rho & = C'^{-1} e^{-\frac{1}{2} \int (|u|^2 + |u_x|^2) dx} \prod_{x \in\T} da(x) db(x).
\end{align}

\smallskip
\noindent {\bf (b)} Construct the measure $\nu$  as a measure absolutely continuous with respect to 
$\rho$ with Radon-Nikodym derivative
\begin{equation}\label{RNungauged}
\frac{d \nu}{d\rho}(u) \,=\, {\mathcal Z}^{-1}   \chi_{\{\|u\|_{L^2} \le B\}}  e^{-\frac{1}{2} N(u)} \,,
\end{equation}
i.e., one needs to show that 
\begin{equation}\label{normalizationungauged}
{\mathcal Z} \,=\, \int  \chi_{\{\|u\|_{L^2} \le B\}}  e^{-\frac{1}{2} N(u)}  d\rho \, < \infty \,.
\end{equation}

Part (b) goes back to the works of Lebowitz, Rose and Speer \cite{LRS} and of Bourgain \cite{B1} for the term 
$\int |u|^6 dx $ part while the integrability of  the term involving   $\int u^2 \bar{u} \bar{u}_x \, dx$ -and hence the construction of $\nu$- is  proved 
in \cite{TTzv};  see also Section 5 in \cite{NORBS}.   Both terms are critical in the sense that integrability 
requires  that $B$ does not exceed  a certain critical value $B^*$. 

Part (a) is a standard problem in Gaussian measures,  treated for example by Kuo \cite{KUO} and 
Gross \cite{GROSS};  see also in \cite{NORBS} for details. Indeed the measure $\rho$ can be realized as an
honest countably  additive Gaussian measure on various Hilbert or Banach spaces  depending on 
one's particular needs.   For example one can construct $\rho$ as the weak limit
of the finite-dimensional Gaussian measures 
\begin{align}\label{finitedim01} 
 d \rho_N & = {\mathcal Z}_{0, N}^{-1} \exp \Big( -\frac{1}{2} \sum_{|n| \leq N} (1+|n|^2)|\ft{u}_n|^2 \Big) 
 \prod_{|n| \leq N} d \ft{a}_n d \ft{b}_n \,,
\end{align}
where $\ft{u}_n = \ft{a}_n + i \ft{b}_n$ is the Fourier transform of $u$.  
For analytical estimates it was convenient in \cite {TTzv} and \cite{NORBS} to consider 
$\rho$  as  measure either on the Hilbert space $H^{\sigma}(\mathbb{T})$ (Sobolev space) for arbitrary 
$\sigma < 1/2$  or as a measure on the Banach space ${\mathcal F} L^{s, r}(\T)$ (Fourier-Lebesgue space 
\cite{Hormander, Gr, Christ}) with norm  $\Vert u \Vert_{ {\mathcal F} L^{s, r}(\T) } \, :=   \Vert \, \langle n  \rangle^s \, 
\widehat u\, \Vert_{\ell^r_n(\Z)}$  and with the conditions $2 \le r <  \infty$ and  $(s-1)r < -1$. 
Note ${\mathcal F} L^{s, r}(\T)$  scales like $H^{\sigma}(\T)$ where $\sigma =  s + \frac{1}{r}- \frac{1}{2} $
and the condition $(s-1)r < -1$ is equivalent to $\sigma<1/2$.  

From a probabilistic stand point, however, and to connect the measure  $\rho$ with the results in 
subsequent sections,  it is also natural to realize this measure on the space of complex-valued 
$2\pi$-periodic continuous   functions $C(\mathbb{T}, \mathbb{C})$. The 
measure  $\rho$ is closely related to the  (complex)   Brownian bridges $Z_{u_o}(x)$ where 
$0\le x \le 2\pi$ and $Z_{u_o}(0)=Z_{u_o}(2\pi)=u_o$.    
Indeed let $\rho( \cdot | u_o)$ denote the measure $\rho$ conditioned on the event $\{u(0)=u(2 \pi)= u_o\}$. 
If $\kappa$  denotes the distribution of $u_o$, then $\kappa$ is a complex Gaussian probability measure
and we have $\rho( \cdot) = \int_{\mathbb{C}}  \rho( \cdot | u_o) d\kappa(u_o)$.  
Then  $\rho( \cdot | u_o)$ is absolutely continuous with respect to the probability distribution $P_{u_o}$ 
of the complex Brownian bridge with 
\begin{equation}\label{rhoBB}
\frac{ d\rho(\cdot | u_o)}{dP_{ u_o}} \,=\,  {\mathcal Z}^{-1}_{u_o} e^{ - \frac{1}{2} \int_0^{2\pi} |u|^2 \, dx} \,.
\end{equation}
This can be easily seen for example by considering the finite-dimensional distribution of $\rho$. 

By combining the results obtained  by Nahmod, Oh, Rey-Bellet, 
and Stafillani in \cite{NORBS}  for the {\em gauged} equation and the results in the present paper we will prove the following theorem:
\begin{theorem}\label{maintheorem1}  The DNLS equation \eqref{DNLS} is $\nu$-almost surely 
well-posed and the measure   $\nu$ is invariant for the flow map $\Psi(t)$ for \eqref{DNLS}. 
\end{theorem}

We now explain why in order to prove Theorem \ref{maintheorem1} one needs to introduce a gauge transformation. We go back to  the existence of (local) solutions to \eqref{DNLS}. By examining the equation one sees there  is  a derivative loss arising from the nonlinear 
term $ (|u|^2 u)_x \, = \,  u^2\,  \overline{u}_x \, +\,  2\, { |u|^2 \, u_x } $  and hence for  low regularity data 
one must somehow make up for this loss.  Since the worse resonant interactions occur on the second term 
${ |u|^2 \, u_x }$ a key idea is to suitably gauge transform the equation to get rid of it,  see 
\cite{HayOz1, HayOz2, Tak1, Herr, GH}.    In the periodic context  a suitable  gauge transformation 
was introduced  by Herr \cite{Herr}.  For $f \in L^2( \mathbb{T})$ let us define 
\begin{equation}\label{thegauge}
G(f)= e^{-iJ(f)}f \,,  \quad  {\rm with} \quad J(f)(x) \,=\, \frac{1}{2\pi}\int_0^{2\pi}\int_\theta^x (|f(y)|^2-m(f))dy\,
d\theta \,,
\end{equation}
and note that  the inverse of $G$ is simply   given by  $G^{-1}(f) = e^{iJ(f)}f$.  Under the gauge $G$, if 
$u$ is a solution  of the DNLS equation \eqref{DNLS} then $w(x,t)=G(u(x,t))$ is a solution  to what we call the GDNLS equation 
\begin{equation}\label{GDNLS+} w_t - i w_{xx} - 2 m(w) w_x = 
- w^2 {\cj w}_x + \frac{i}{2} |w|^4 w  - i  \psi(w) w - i m(w) |w|^2 w  \,,
\end{equation} 
where
\begin{equation*}
\psi(w):= -\frac{1}{\pi} \int_{\T}  {\rm Im}  ( w {\cj w}_x ) \, dx \, + \frac{1}{4 \pi} \int_{\T} | w|^4 dx - m(w)^2  \,.
\end{equation*}

The main result of the present paper is to show how the measure $\nu$ is transformed under the gauge 
transformation $G$. The image of $\nu$ under $G$ is denoted \footnote[6]{This $\mu$ is not yet the same as the $
\mu$ constructed in \cite{NORBS} that we referred to in the Introduction. But it will indeed be the same as a 
consequence of \eqref{RNgauged} after we prove Theorem \ref{maintheorem2}.} 
by $\mu$ and is, by definition, given  by 
\begin{equation}\label{mudef}
\mu(A) \,:= \, \nu( G^{-1}(A) ) \,=\, \nu\left( \{x ; G(x) \in A\}\right),
\end{equation}
where $A$ is any measurable set.  We will use the notation $\mu =\nu \circ G^{-1}$ in the sequel. 
We have

\begin{theorem}\label{maintheorem2}  For sufficiently small $B$, the measure $\mu=\nu\circ G^{-1}$ is 
absolutely continuous with respect to the Gaussian measure $\rho$ and we have 
\begin{equation} \label{RNgauged}
\frac{d \mu}{d \rho}(w) =  \tilde{\mathcal Z}^{-1}  \chi_{\{\|w\|_{L^2} \le B \}}\, e^{- \frac{1}{2} \mathcal{N}(w)} \,,
\end{equation}
where 
$$
\mathcal{N} (w) \,=\,  - \frac{1}{2} {\rm Im} \int w^2 \cj{w} \cj{w_x} \,dx + 2m(w)\, {\rm Im}\, \int w \cj{w_x} \, dx  
 -\frac{1}{2}m(w)\,  \int |w|^4 \, dx  +   2\pi m(w)^3,
$$
and $ \tilde{\mathcal Z}$ is a normalization constant.  
\end{theorem}

For the measure $\mu$, as given by \eqref{RNgauged},  the following result is proved in \cite{NORBS}, 
see Theorem 6.3,  6.5, 7.1, and 7.2.

\begin{theorem}\label{maintheorem3}  
The GDNLS equation \eqref{GDNLS+} is $\mu$-almost surely well-posed and the measure  
$\mu$ is invariant for the flow map $\Phi(t)$ for \eqref{GDNLS+}. 
\end{theorem}

\begin{remark}{\rm 
In \cite{Herr} and \cite{NORBS} one actually performs another supplementary transformation 
to get rid of the term  $2 m(w) w_x $ on the left hand side of \eqref{GDNLS+}. Indeed 
if we set $v(x,t ) =  w( x - 2 tm(w), t)$ then $v$ is a solution of 
\begin{equation}\label{GDNLS} v_t - i v_{xx}  = 
- v^2 {\cj v}_x + \frac{i}{2} |v|^4 v  - i  \psi(v) v - i m(v) |v|^2 v  \,.
\end{equation} 
A simple argument given in section 7 of \cite{NORBS} show that the measure $\mu$ is invariant
for both the flow maps for \eqref{GDNLS+} and \eqref{GDNLS}.  
}
\end{remark}

To conclude one notes that Theorem \ref{maintheorem1} follows immediately from 
Theorem \ref{maintheorem2} and \ref{maintheorem3}.  We are thus left to prove Theorem  \ref{maintheorem2}.

\smallskip 

\subsection{An heuristic introduction of $\mu$} \label{heuristic} To understand the form of the measure $\mu$ 
we give here a purely heuristic argument,  a rigorous proof 
is given in the next section. Let us first recall how the invariants for  DNLS transform under $G$. 
Since  $u = e^{iJ(w)} w$  we have  $m(u)=m(w)$ and   $u_x = e^{iJ(w)}( w_x +  i J(w)_x w)$ with 
$J(w)_x =  |w|^2 - m(w)$ and we  obtain after straightforward computations
\begin{eqnarray}
H(u)  &=& {\rm Im} \int_{\T} u {\cj u}_x \, dx  + \frac{1}{2} \int_{\T} |u|^4 \, dx
\nonumber \\
&=&  {\rm Im} \int_{\T} w  {\cj w}_x   - \frac{1}{2} \int_{\T} |w|^4 \, dx  + 2\pi m(w)^2=:\mathscr H(w) \,, \label{oo4}
\end{eqnarray}
and
\begin{eqnarray}
u_x \cj{u_x}  \,&=&\,   \left( {w_x} + i J(w)_x w \right)    \left( \cj{w_x} - i J(w)_x \cj{w} \right)  \nonumber  \\
\,&=&\,  w_x \cj{w_x}   -2  {\rm Im} \, w^2  \cj{w} \cj{w_x}   +2m(w) \, {\rm Im} \,w \cj{w_x}   
            + \left( |w|^6 - 2m |w|^4 +  m(w)^2\,|w|^2\right)      \,,                
                          \label{oo2}
                       \end{eqnarray}
as well as
\begin{equation}\label{oo3}
u^2\cj{u} \cj{u_x} =w^2\cj{w} \cj{w_x} - i|w|^6 +i m(w)\, |w|^4 \,.
\end{equation}                      
Hence by using  \eqref{energy}, \eqref{oo2},  \eqref{oo3}  we find  
\begin{eqnarray} \label{equiv-energy}
E(u) \,&=&\,  \int w_x \cj{w_x} \, dx  - \frac{1}{2} {\rm Im} \int w^2 \cj{w} \cj{w_x} \,dx 
+ 2m(w) \, {\rm Im}\, \int w \cj{w_x} \, dx  \\
&&  -\frac{1}{2}m(w)\,  \int |w|^4 \, dx  +  2\pi m(w)^3 \nonumber \\
\,&=:&\, \EE(w) \,. 
\end{eqnarray}

\begin{remark}{\rm  
Notice that $\EE(w)$ involves the other conserved quantities $\mathscr H(w)$ and $m(w)$ and if we define 
\begin{equation}\label{ew}
\mathscr E(w):=  \int_{\T} | w_x|^2 \, dx - \frac{1}{2} {\rm Im} \int_{\T} w^2 {\cj {w} \, \cj{w_x}} \, dx 
+ \frac{1}{4 \pi}  \biggl(\int_{\T} |w(t)|^2 \, dx \biggr)\biggl(\int_{\T} |w(t)|^4 \, dx\biggr),
\end{equation}
we then have 
\begin{equation}\label{eu}
E(u) = \mathscr E(w) + 2m\, \mathscr H(w) - 2\pi \, m^3=\EE(w),
\end{equation}
and so $\mathscr E(w)$ is also a conserved quantity. One could build invariant measures using 
$\mathscr E(w)$ rather than $\EE(w)$ but they would turn out to be equivalent measures.
}
\end{remark} 

\medskip 

Let us pretend that the measure $\nu$ is the  measure  with density  
\begin{equation}
\chi_{\{\|u\|_{L^2} \le B\}}  e^{-\frac{1}{2} N(u)} e^{-\frac{1}{2} \int (|u|^2 + |u_x|^2) dx}  
= \chi_{\{\|u\|_{L^2} \le B\}}  e^{-\frac{1}{2} E(u) -\frac{1}{2} \int |u|^2 dx} 
\end{equation} 
with respect to the (nonexistent!)  infinite dimensional  Lebesgue measure $\prod_{x \in\T} da(x) db(x)$.  
Let us  assume furthermore that this nonexistent Lebesgue measure is left invariant under $G$.  
Then we would simply obtain from \eqref{equiv-energy} that 
\begin{eqnarray} 
d \mu \,&=&\, {\tilde C}^{-1} \chi_{\{\|w\|_{L^2}\le B\}} \, e^{- \frac{1}{2} \EE(w)  -\frac{1}{2} \int |w|^2 dx} \, 
\prod_{x \in \mathbb{T}} d  a(x) d b(x) \\
&=&\, {\mathcal Z}^{-1} \chi_{\{\|w\|_{L^2}\le B\}} \, e^{- \frac{1}{2} {\mathcal N}(w)} d \rho \,, \nonumber
\end{eqnarray} where 
$$
\mathcal{N} (w) \,=\,  - \frac{1}{2} {\rm Im} \int w^2 \cj{w} \cj{w_x} \,dx 
+ 2m(w)\, {\rm Im}\, \int w \cj{w_x} \, dx   -\frac{1}{2}m(w)\,  \int |w|^4 \, dx  +   2\pi m(w)^3 
$$
is the nonquadratic part of the energy $\EE(w)$.  

The crucial problem to understand rigorously the transformation of $\mu$ under $G$ is actually 
to understand  the transformation of the Gaussian part $\rho$ of $\mu$ under $G$ since the transformation of the 
weight is computed easily as in the formal computation above.  This is achieved in Theorem \ref{T1} below where  
in order to analyze the transformation of $\rho$ the main ingredients will be:

\smallskip 
\noindent {(i)} The relation  between $\rho$ and Brownian bridges,  see eq. \eqref{rhoBB}.  

\smallskip
\noindent {(ii)} The well-known fact that a Brownian bridge can be obtained by conditioning a Brownian motion 
to return at  its starting point. 

\smallskip
\noindent {(iii)} The conformal invariance of Brownian motion. Note that since $w = e^{-iJ(u)} u,$ 
$J(u)=J(w)$ and  $J(u)=J(|u|),$ it is more convenient to consider this transformation in terms of 
the variables  
\begin{equation}\label{linearization}
|w| = |u|   \,, \qquad   \arg(w) = \arg(u) -i J(|u|).
\end{equation}
By conformal invariance of Brownian motion and (ii),  $|u|$ and $\arg(u)$ have a Gaussian distribution after 
a suitable reparametrization.  The transformation \eqref{linearization} is easy to understand.  In particular
if we condition on $| u | $,  the transformation of $\arg(u)$ is a simple translation by a fixed vector which leads to 
the next item (iv).

\smallskip
\noindent{(iv)} The Cameron-Martin formula for the transformation of Gaussian measure under a translation by 
a fixed vector,  see e.g.  \cite{CM,Bogachev}. We will use here the following special case of the Cameron-Martin
theorem.

\begin{theorem}\label{CMBB}{(Cameron-Martin theorem for real Brownian bridge)}.   Let $X(s)$, $0 \le s \le S$, 
be a real  Brownian bridge, with $X(0)=X(S)=x_o$ and law $R_{x_o}$.  Let $k(s)$, $0 \le s \le S$, 
be an  absolutely continuous real-valued function such that $k(0)=k(S)=k_o$ and 
$\int_0^S |k'(s)|^2 \, ds < \infty$.  Then the law  of $X(s) + k(s)$ is absolutely continuous with respect to 
$R_{x_o + k_o}$ with Radon-Nikodym derivative
$$
\exp \left(  \int_0^S k'(s) dX(s) - \frac{1}{2} \int_0^S |k'(s)|^2 ds \right) \,.
$$
\end{theorem}

\section{Brownian bridges under gauge transformations}\label{general} 

Let $Z_{u_o}(x)$ be a standard complex Brownian bridge{\footnote[7]{Although in the probability literature 
it is customary to use $t$ as the variable of Brownian motions, here we use $x$ instead in order not to confuse it with the 
time $t$ of the evolution equation \eqref{DNLS}.}}  on the interval $0 \le x \le 2 \pi$ and with 
$Z_{u_o}(0)=Z_{u_o}(2\pi) = u_o$. The law of $Z_{u_o}$ is denoted by $P_{u_o}$ and is a 
Gaussian probability on $\{ Z \in C(\mathbb{T}\,;\, \mathbb{C}), Z(0)=u_o\}$.  Since no confusion arises we 
will omit the index $u_o$ in the sequel and denote the Brownian bridge simply  by  $Z$ and 
its probability distribution by $P$.   We consider first the transformation of a  complex Brownian bridge 
under a class of  transformations which contains in particular  the gauge transformation $G$ given 
in \eqref{thegauge}.

We assume that the map $G$ satisfies the following condition

\smallskip
\noindent 
{\bf (C)} The map 
$ G \,:\, C( \mathbb{T}, \mathbb{C}) \to C( \mathbb{T}, \mathbb{C}) $
has the form 
\begin{equation}\label{gauge}
G(Z)(x) \,=\,   e^{- i J(Z)(x)} Z(x)  \,,
\end{equation}
where  $J:  C( \mathbb{T}, \C)  \to C(\mathbb{T}, \R)$ depends only on $|Z|$ and is such that 
\begin{equation}\label{ju}
\frac{d}{dx}J(Z)(x)=h(|Z|)(x)\,,
\end{equation}
and $h(|Z|)(x)$ is continuous in $x$ for $P$-almost all choices of the process $|Z|$.

\smallskip The gauge transformation \eqref{thegauge} in Section \ref{intro} satisfies condition {\bf (C)} since 
we have
\begin{eqnarray} 
J(Z)(x) \,&=&\, \frac{1}{2\pi}\int_0^{2\pi}\int_\theta^x \left(|Z(y)|^2- \frac{1}{2\pi} \int_0^{2\pi} |Z(\xi)|^2 \, d\xi \right) \,dy \,
d\theta  \,, \\
h(|Z|)(x) \,&=&\,  |Z(x)|^2- \frac{1}{2\pi} \int_0^{2\pi} |Z(\xi)|^2 \, d\xi\,.
\end{eqnarray}

\begin{theorem}\label{T1} Let $Z$ be a standard complex  Brownian bridge 
with law $P$.  Let $G$ be a map which satisfies the condition {\bf (C)}.  
If we have 
\begin{equation}\label{novikov}
\E_{P} \left[\exp\left(\mbox{\rm Im}\int _0^{2\pi} h(|Z|)(x) Z(x) \,d\cj{Z}(x) -\frac{1}{2} \int _0^{2 \pi}|h(|Z|)(x)|^2|Z(x)|^2\,dx \right)\right] = 1 \,,
\end{equation}
then   $P \circ G^{-1}$ is absolutely  continuous with respect to the  Brownian bridge $P$ with 
Radon-Nikodym derivative
\begin{equation}\label{RNbridge}
\exp\left(\mbox{\rm Im}\int _0^{2\pi} h(|Z|)(x) Z(x) \,d\cj{Z}(x) -\frac{1}{2} \int _0^{2 \pi}|h(|Z|)(x)|^2|Z(x)|^2\,dx \right) \,.
\end{equation}
\end{theorem}

\begin{remark}{\rm  Assume for a moment that we are not in a periodic setting, that  $Z$ is a standard 
Brownian motion instead of a Brownian bridge, and that  
\begin{equation}\label{non-anti}J(Z)(x) = \int_0^x f(|Z|(t)) \, dt \,,
\end{equation}
for some  real-valued 
continuous function  $f $.  The expression in the right hand side 
of \eqref{non-anti} looks very similar to the corresponding one in \eqref{ju}, but it is actually easier to handle since it is non anticipative, in the sense that it depends on the Brownian motion up to ``time" $x$ and not later.
Thanks to this fact  the  Radon-Nikodym  derivative for the transformation $G$ can be 
computed as a consequence of Girsanov  formula (e.g. \cite{Oks}).   Indeed if we set  
$R= J(Z)$ then we have  
$$
\tilde Z \equiv G(Z)=  e^{-iR} Z \,.
$$
By Ito's formula we have $dR\,=\,  f( |Z|) dx$ and 
\begin{eqnarray}
d\tilde Z  \,&=&\, -i e^{-iR} Z dR  + e^{-iR} dZ   - \frac{1}{2} e^{-iR}Z  dR^2  + i e^{-iR} dR dZ + 0 \frac{1}{2} dZ^2 \nonumber  \\ 
\,&=&\,  e^{-iJ(Z)}(- i Z) f( |Z|) dx  + e^{-iJ(Z)} dZ \nonumber \\ 
\,&=&\,  - i {\tilde Z}  f( |{\tilde Z}|) dx  + e^{-iJ(Z)} dZ \,, \nonumber
\end{eqnarray}
where we have used that $|\tilde Z| = |Z|$.  Since $J(Z)(x)$ is a nonanticipating functional   
$e^{-iJ(Z)} dZ$ is a Brownian motion (see \cite{Oks})  and therefore 
${\tilde Z}$ has the same law as the  solution of the SDE
$$
d {\tilde Z}  \,=\,   -i {\tilde Z}  f( {\tilde Z}) dx  + dZ \,.
$$
An application of Girsanov Theorem gives now the form of the Radon-Nikodym derivative as in 
\eqref{RNbridge}. 
}
\end{remark}

\medskip
\noindent
{\em Proof of Theorem \ref{T1}} The remark above explains why the kind of gauge transformations we consider cannot be studied  directly by the Girsanov Theorem and some manipulation needs to be performed.

In the course of the proof we will use some properties of the complex Brownian motion which we denote by 
$B(x)$  (again we omit from the notation  the choice of $B(0)$.)  We recall first the  well-known fact, see e.g. 
\cite{Oks}, that a Brownian bridge $Z(x)$ is obtained from a  Brownian  motion by conditioning $B$ on 
the event $\{ B(2\pi)=B(0)\}$.    

Furthermore  we will use the conformal invariance of Brownian motion,  that is if  $A= A_1+ i A_2$ is 
a complex Brownian motion, and $\phi$  is analytic function  then $B= \phi(A)$ is, after a suitable time 
change, again a complex  Brownian motion (see e.g. \cite{Oks}, Example 8.22).  For $B(x) = \exp(A(s))$ 
the  time change is  given by 
\begin{equation}\label{ts} x=x(s)=\int_0^s \left|e^{A(r)}\right|^2 dr \,=\,\int_0^s e^{2 A_1(r)} dr  \qquad \frac{dx}{ds}=\left|e^{A(s)}\right|^2=
|B(x(s))|^2,
\end{equation}
or equivalently 
\begin{equation}\label{st}s(x)=\int_0^x \frac{dr}{|B(r)|^2}, \qquad \frac{ds}{dx}=\frac{1}{|B(x)|^2} \,.
\end{equation}
If we write  $B(x)$ in polar coordinate 
\begin{equation}
B(x) = |B(x)| e^{i \arg(B)(x)} \,,
\end{equation}
we have 
\begin{equation}
A(s) =  A_1(s) + i A_2(s) = \log |B(x(s))| + i \arg(B) (x(s)) \,,
\end{equation}
and $A_1$ and $A_2$ are real independent Brownian motions.  

In view of conditioning the Brownian motion to obtain a Brownian bridge we are interested in 
$B(x)$  for  $0 \le x \le 2\pi$ and  thus we introduce the  stopping time
\begin{equation}
S= \inf \left\{ s \,;\, \int_0^{s} \left| e^{A(r)}\right|^2 dr=2 \pi \right\} \,,
\end{equation}
and remark,  for future use, the important fact that the stopping time $S$ depends only 
on the real part $A_1(s)$ of $A(s)$, or equivalently only $|B|(x).$ 

For the Brownian bridge $Z(x)$ let us define $W(s)=W_1(s)+i W_2(s)$ where
\begin{equation}
W_1(s) = \log |Z(x(s))| \,, \quad  W_2(s) = \arg (Z) (x(s)) \,,
\end{equation}
where $x(s)=\int_0^s e^{ 2W_1(r)} dr$ and $0 \le x \le 2\pi$ is equivalent to $ 0 \le s \le S$.  Furthermore we denote by $Q$ the law of $W$.
The stopping time depends only $|Z|$, so once we condition on the process $W_1(s)=\log|Z(x(s))|$,   
the conditional law of the process $W_2(s)=\arg(Z) (x(s))$  is now the law of   a real 
Brownian bridge on the interval  $0 \le s \le S$.   

If we define ${\widetilde W}={\mathcal L(W)}={\widetilde W}_1 + i {\widetilde W}_2 $ by 
\begin{eqnarray}
{\mathcal L}(W)(s)&:=&  W_{1}(s) + i \left[ W_{2}(s) - J(|Z|)(x(s)) \right] \nonumber  \\
&=&  W_1(s) + i \left[ W_2(s) - J (e^{W_1})(s) \right]  \label{ltransformation} \,,
\end{eqnarray}
then we have 
$$
e^{ {\mathcal L}(W)(s)} = G(Z)(x(s)),
$$ 
where $G$ is as in \eqref{gauge}. 
From  \eqref{ltransformation} we observe first that the (marginal) law of $W_1$, under $Q$, and 
the (marginal) law of the real part  ${\widetilde W}_1$ of ${\mathcal L}(W)$,  under 
$Q\circ {\mathcal L}^{-1}$, coincide.   Furthermore the conditional law of $W_2$, conditioned on 
the real part $W_1$, is the law of  a real Brownian bridge  on the interval $[0,S]$. 
Therefore from \eqref{ltransformation}   the conditional law of ${\widetilde W}_2$, conditioned on 
the real part $W_1={\widetilde W}_1$,  can be computed   using Cameron-Martin 
formula, Theorem \ref{CMBB},  since ${\widetilde W}_2$ is obtained from $W_2$ by translating by a function    which  depends only on $W_1$.


Since 
\begin{equation}
\frac{dJ(u(x(s)))}{ds}=(J(u))'(x(s))\frac{dx(s)}{ds}=h(|u|)(x(s))\frac{dx(s)}{ds} \,,
\end{equation}
the Cameron-Martin formula implies that the Radon-Nikodym derivative of the law of imaginary part 
of ${\mathcal L}(W)$  with respect to the law of a real brownian bridge on the interval 
$0 \le s \le S$ is given by 
$$
\exp \left(  \int_0^S   h( |Z|) (x(s))  \frac{dx}{ds}  dW_2(s) 
-\frac{1}{2} \int_0^S  h( |Z|)^2 (x(s)) \left(\frac{dx}{ds}\right)^2 ds    \right)  \,.
$$

To conclude we finally  re-express the Radon-Nikodym derivative in terms of $Z$ and $x$. We have 
\begin{equation}
dW(s)=\frac{1}{Z(x(s))} \frac{dx}{ds} d Z (x(s)) =\cj Z (x(s)) dZ (x(s)) \,,
\end{equation}
and thus  
\begin{equation}\label{term1}
 \int_0^S   h( |Z|) (x(s))  \frac{dx}{ds}  dW_2(s) \,=\, {\rm Im} \int_0^{2 \pi} h(|Z|)(x) \cj{Z} \, dZ(x) \,,
\end{equation}
and 
\begin{equation}\label{term2} \int_0^S  h( |Z|)^2 (x(s))  \left( \frac{dx}{ds} \right)^2  ds   =\int _0^{2 \pi}  |h(|Z|)(x)|^2 |Z(x)|^2\,dx.
\end{equation}
This concludes the proof of Theorem \ref{T1}.  \qed

%

\subsection{Application to the periodic derivative NLS} 

Let us consider the measure $\mu$ given in the introduction, see \eqref{WienerUG} and 
\eqref{RNungauged}.   In this section we prove Theorem \ref{maintheorem2}  using Theorem \ref{T1}.

\noindent 
{\em Proof of Theorem \ref{maintheorem2}}  
We note that, using \eqref{oo3}, the Radon-Nykodym derivative 
$\frac{d \nu}{d \rho}$ transforms under the  gauge  $G$ as  
\begin{equation}\label{densitytransform}
N(G^{-1}(w)) = \frac{3}{2} {\rm Im} \int w^2\cj{w} \cj{w_x} \, dx  - \int |w|^6\, dx + \frac{3}{2} m \int |w|^4 \, dx \,.
\end{equation}
Furthermore by the results of \cite{B1,TTzv,NORBS} $ \chi_{\{\| w\|_{L^2} \le B\}} e^{ -\frac{1}{2} N(G^{-1}(w))}  
\in L^1( \rho)$ for sufficiently small $B$. Therefore it is enough to consider how $\rho$ transforms under the gauge transformation, i.e., we consider 
the measure $\tilde{\rho} = \rho \circ G^{-1}$.  Without cutoff on $\|u\|_{L^2}$ one cannot expect $\tilde{\rho}$ 
to be absolutely continuous with respect to  $\rho$,  but all we really need is that the restriction of $\tilde{\rho}$ on 
$\{m(w) \le \frac{B^2}{2\pi} \}$  be  absolutely continuous with respect to  $\rho$. Alternatively we  can  
incorporate the cutoff in $J$ by redefining $J$ to be  
\begin{equation}\label{judnls}
J(w)(x):=\left\{
\begin{array}{ll}
 \frac{1}{2\pi}\int_0^{2\pi}\int_\theta^x (|w(y)|^2-m(w))dy\,d\theta &\qquad \mbox{ if } m(w) \leq \frac{B^2}{2\pi}\\
\qquad\\
0 &\qquad \mbox{otherwise}.
\end{array}\right.
\end{equation}
so that  we have 
\begin{equation}\label{j'}
\frac{d}{dx}J(w)(x)=h(|w|)(x)=\left\{
\begin{array}{ll}
 |w(x)|^2-m(w) &\qquad \mbox{ if } m(w) \leq \frac{B^2}{2\pi} \\
\qquad\\
0 &\qquad \mbox{otherwise}.
\end{array}\right.
\end{equation}
By the results in \cite{B1,TTzv,NORBS}  
\begin{eqnarray}
&&\exp\left( {\rm Im} \int h(|w|) w \cj{w}_x \, dx - \frac{1}{2} \int h(|w|)^2 |w|^2  \, dx \right) =\label{rhononlin} \\
&&  \chi_{\{ m(w) \le \frac{B^2}{2\pi} \} }  \exp\left(  {\rm Im} \int (|w|^2 - m(w))  w  \cj{w}_x \, dx  -\frac{1}{2} \int ( |w(x)|^2- m(w))^2 |w(x)|^2\,dx
\right)  \nonumber
\end{eqnarray}
belongs to $L^1(\rho)$ for sufficiently small $B$.     By conditioning on $\{u(0)= u(2 \pi)= u_o\}$ 
and using equation \eqref{rhoBB}  we conclude that the Novikov condition  \eqref{novikov} is satisfied for 
almost every $u_o$.   Therefore using  Theorem \ref{T1}  the Radon-Nikodym derivative 
$\frac{ d\tilde{\rho}}{ d\rho}$  is given by $\eqref{rhononlin}$. Finally combining  the equations 
\eqref{rhononlin}  and \eqref{densitytransform} we  obtain equation \eqref{RNgauged}.  
\qed

\bigskip

{\bf Acknowledgment.}  Andrea R. Nahmod and Gigliola Staffilani would like to warmly thank 
the Radcliffe Institute for Advanced Study at Harvard University for its wonderful hospitality
while part of this work was being carried out. They also thank their fellow Fellows for the 
stimulating environment they created.

\end{document}